\documentclass[12pt]{amsart}

\usepackage{amsmath,amsthm,amssymb}
\font\teneufm=eufm10
\font\seveneufm=eufm7
\font\fiveeufm=eufm5
\newfam\frakturfam

\textfont\frakturfam=\teneufm
\scriptfont\frakturfam=\seveneufm
\scriptscriptfont\frakturfam=\fiveeufm

\newtheorem{lm}{Lemma}
\newtheorem{theor}{Theorem}
\newtheorem{co}{Corollary}

\newtheorem{prob}{Problem}
\def\bee{\begin{eqnarray}}
\def\bes{\begin{eqnarray*}}
\def\eee{\end{eqnarray}}
\def\ees{\end{eqnarray*}}

\def\a{\alpha}

\def\Proof{{\sl Proof.}\ }

\title{Polarization algebras and their relations}

\begin{document}
\date{}
\maketitle

\begin{center}

{\bf Ualbai Umirbaev}\footnote{Supported by an NSF grant DMS-0904713 and by an MES grant 0755/GF of Kazakhstan; Eurasian National University,
 Astana, Kazakhstan and
 Wayne State University,
Detroit, MI 48202, USA,
e-mail: {\em umirbaev@math.wayne.edu}}

\end{center}

\begin{abstract} Using an approach to the Jacobian Conjecture by L.M. Dru\.zkowski and K. Rusek \cite{DR85}, G. Gorni and G. Zampieri \cite{GZ96}, and A.V. Yagzhev \cite{Yagzhev00-1}, we describe a correspondence between finite dimensional symmetric algebras and homogeneous tuples of elements of polynomial algebras. We show that this correspondence closely relates Albert's problem \cite[Problem 1.1]{Dnestr} in classical ring theory and the homogeneous dependence problem \cite[page 145, Problem 7.1.5]{Essen}  in affine algebraic geometry related to the Jacobian Conjecture. We demonstrate these relations in concrete examples and formulate some open questions.
\end{abstract}


\noindent


\section{Introduction}

\hspace*{\parindent}

The main objective of this paper is to connect two groups of specialists who are working on a closely connected problems and sometimes have intersections. Joining the efforts of these groups may be fruitful in studying the Jacobian Conjecture \cite{Essen}. Namely,  Albert's problem \cite[Problem 1.1]{Dnestr} in classical ring theory and the homogeneous dependence problem \cite[page 145, Problem 7.1.5]{Essen}  in affine algebraic geometry are closely connected. Ring theory specialists are focused in studying only binary algebras and have some positive results in small dimensions. Affine algebraic geometry specialists have also some positive results in small dimensions and have some negative results for $m$-ary algebras if $m\geq 3$.

To make the subject more intriguing let's start with two examples.
A well known quadratic homogeneous automorphism \cite[page 98]{Essen}
\bee\label{Zf1}
(x-ys,y+zs-xt,z-yt,s,t)
\eee
of the affine space $\mathbb{A}^5$ over a field $k$ corresponds to the Suttles'  example of a commutative power-associative $5$ dimensional algebra which is solvable but not nilpotent \cite{Suttles}.  Another well known  homogeneous automorphism \cite[page 97]{Essen}
\bee\label{Zf2}
(x+s(xt-ys),y+t(xt-ys),s,t)
\eee
of the affine space $\mathbb{A}^4$ gives an example of a ternary symmetric power-associative $4$ dimensional algebra which is also solvable but not nilpotent. Recall that (\ref{Zf1}) is a tame automorphism and (\ref{Zf2}) is the first candidate to be non-tame in the case of $4$ variables.

The paper is organized as follows. In Section 2 we give the definition of polarization algebras and prove some elementary results. In Section 3 we define symmetric algebras with Engel condition and describe  an approach to the Jacobian Conjecture by L.M. Dru\.zkowski and K. Rusek \cite{DR85}, G. Gorni and G. Zampieri \cite{GZ96}, and A.V. Yagzhev \cite{Yagzhev00-1}. In Section 4 we describe a connection between Albert's problem and the homogeneous dependence problem and consider some examples. In Section 5 we formulate some open problems.

\section{Polarization algebras}

\hspace*{\parindent}

Let $k$ be an arbitrary field of characteristic zero.
It is reasonable to call an algebra $A$ over $k$ with one $m$-ary multilinear operation $\langle \cdot,\cdot,\ldots, \cdot\rangle$ {\em  symmetric} if
\bes
\langle x_1,x_2,\ldots,x_m\rangle=\langle x_{\sigma(1)},x_{\sigma(2)},\ldots,x_{\sigma(m)}\rangle
\ees
for any $x_1,x_2,\ldots,x_m\in A$ and for any $\sigma\in S_m$, where $S_m$ is the symmetric group on $m$ symbols.
If $m=2$ then a symmetric algebra becomes a commutative (non-associative) algebra.

Let $P_n=k[x_1,x_2,\ldots,x_n]$ be the polynomial algebra over $k$ in the variables  $x_1,x_2,\ldots,x_n$. Denote by $P_n^n$ the set of columns over $P_n$ of height $n$. For convenience we write columns as rows.
Put $V=k^n$. Every $n$-tuple $H=(h_1,h_2,\ldots,h_n)\in P_n^n$ of elements of  $P_n$ uniquely defines the polynomial mapping $H : V \rightarrow V$. We often say that $H$ is a polynomial mapping or endomorphism without mentioning $V$. The $n$-tuple $X = (x_1,x_2,\ldots,x_n)$ represents the identity mapping. We also use notations $P_n=k[X]$, $f=f(X)\in P_n$, and $H=H(X)$.

For any $n$-tuple $H=(h_1,h_2,\ldots,h_n)$ we put
\bes
\deg(H)= \max_i(\deg(h_i)),
\ees
where $\deg$ is the standard degree function on $P_n$.

 Let $H=H(X)=(h_1,h_2\ldots,h_n)$ be an $n$-tuple of homogeneous elements of $P_n$ of degree $m\geq 2$, i.e., all $h_i$ are homogeneous of degree $m$. Notice that $H$ is written in the form of a row vector but we always consider the vectors as column vectors. Consider the $n$ dimensional vector space $V=k^n$.  Let $\{x_{ij} |1\leq i\leq m,  1\leq j\leq n\}$ be a set of independent commutative variables. Put $X_i=(x_{i1},x_{i2},\ldots,x_{in})$ for all $i$.
 The polarization \cite{Dolgachev,Kraft2} (or full linearization)
\bes
\langle X_1,X_2,\ldots,X_m\rangle= \frac{1}{m!} \sum_{\mu\subseteq [m]} (-1)^{m-|\mu|} H(\sum_{i\in \mu}X_i)
\ees
of $H$, where $[m]=\{1,2,\ldots,m\}$ and $|\mu|$ is the cardinality of $\mu$, defines a unique $m$-ary symmetric operation on $k^n$. This algebra is called the {\em polarization} algebra of $H$ and will be denoted by $P_H$.

Obviously, the restitution  $\langle X,X,\ldots, X\rangle$ gives again $H$ \cite{Dolgachev,Kraft2}.

Let $A$ be an arbitrary symmetric $n$ dimensional algebra over $k$ and let
\bes
(f)=\{f_1,f_2,\ldots,f_n\}
\ees
 be a linear basis for $A$. Then
 \bes
  Y=x_1\otimes f_1+x_2\otimes f_2+\ldots+x_n\otimes f_n
 \ees
is called a {\em generic} element of $A$. It is well-known (see, for example \cite{KBKA}) that the $k$-subalgebra generated by $Y$ in $k[X]\otimes_k A$ is a free one generated algebra of the variety of algebras generated by $A$. Consider the element
\bes
\langle Y,Y,\ldots, Y\rangle=g_1\otimes f_1+g_2\otimes f_2+\ldots+g_n\otimes f_n
\ees
where $g_1,g_2\ldots,g_n\in k[x_1,x_2,\ldots,x_n]$. Note that $G=(g_1,g_2,\ldots,g_n)$ is a homogeneous $n$-tuple of  degree $m$. The $n$-tuple
$G$ is called the {\em $n$-tuple of $A$ corresponding to $(f)$}. If the basis $(f)$ is not specified then we say that $G$ is an $n$-tuple of $A$.

\begin{lm}\label{Zl1}  Let $A$ be an arbitrary $m$-ary symmetric $n$ dimensional algebra over $k$ and let $G$ be an $n$-tuple of $A$. Then $A$ is isomorphic to the polarization algebra $P_G$.
\end{lm}
\Proof Let $G$ be the $n$-tuple of $A$ corresponding to a basis $(f)$. Denote by
\bes
(e) = \{e_1,e_2,\ldots,e_n\}
\ees
the standard basis for $V=k^n=P_G$ (columns of the identity matrix of order $n$). Let $\varphi : A\rightarrow P_G$ be the linear mapping defined by $\varphi(f_i)=e_i$ for all $i$.

It is sufficient to show that $\varphi$ is a homomorphism of algebras, i.e.,
\bes
\varphi\langle v_1,v_2,\ldots,v_m\rangle=\langle \varphi(v_1),\varphi(v_2),\ldots,\varphi(v_m)\rangle
\ees
for all $v_1,v_2,\ldots,v_m\in V$. This is equivalent to
\bee\label{Zf3}
\varphi\langle v,v,\ldots,v\rangle=\langle \varphi(v),\varphi(v),\ldots,\varphi(v)\rangle
\eee
for all $v\in V$ since the both operations are  $m$-ary symmetric and $k$ is a field of characteristic zero.

Put $v=\lambda_1 f_1+\lambda_2f_2+\ldots+\lambda_nf_n\in A$ and
$w=\lambda_1 e_1+\lambda_2e_2+\ldots+\lambda_ne_n=(\lambda_1,\lambda_2,\ldots,\lambda_n)\in P_G$.
We get $\varphi(v)=w$ by the definition of $\varphi$. We also have
\bes
\langle v,v,\ldots, v\rangle=g_1(\lambda)\otimes f_1+g_2(\lambda)\otimes f_2+\ldots+g_n(\lambda)\otimes f_n,
\ees
where $g_i(\lambda)=g_i(\lambda_1,\lambda_2,\ldots,\lambda_n)$ for all $i$, since $G$ is the $n$-tuple of $A$ corresponding to $(f)$. By the definition of the product in $P_G$, we get
\bes
\langle v,v,\ldots, v\rangle=g_1(\lambda)\otimes e_1+g_2(\lambda)\otimes e_2+\ldots+g_n(\lambda)\otimes e_n.
\ees
Consequently, (\ref{Zf3}) holds. $\Box$

Thus, every $m$-ary symmetric finite dimensional algebra is a polarization algebra. Now we describe the conditions when two polarization algebras are isomorphic.

Let $\mathrm{Aut}(\mathbb{A}^n)$ be the group of all polynomial automorphisms of the the affine space $\mathbb{A}=k^n=V$. The group $\mathrm{Aut}(\mathbb{A}^n)$ acts on the set of all polynomial mappings on $V$ by conjugation. Denote by $GL_n(k)$ the group of all linear automorphisms of $V$. If $H$ is homogeneous of degree $m$ then $\a H \a^{-1}$ is homogeneous of degree $m$ for all $\a\in GL_n(k)$.
Two homogeneous endomorphisms $H$ and $G$ of the same degree $m$  are called {\em linearly conjugate} if there exists $\a\in GL_n(k)$ such that $\a H \a^{-1}=G$.
\begin{lm}\label{Zl2}  Let $H$ and $G$ be homogeneous $n$-tuples of elements of $P_n$ of the same degree $m\geq 2$. Then
the polarization algebras  $P_H$ and $P_G$ are isomorphic if and only if $H$ and $G$ are linearly conjugate.
\end{lm}
\Proof Let $\mathrm{Hom}(V,V)$ be the set of all linear transformations of $V$. As in the proof of Lemma \ref{Zl1},
 $\varphi\in \mathrm{Hom}(V,V)$ is a homomorphism of algebras $P_H$ and $P_G$ if and only if (\ref{Zf3}) holds.
 Consequently,  $\varphi$ is a homomorphism if and only if $\varphi(H(v))=G(\varphi(v))$, i.e., $\varphi\circ H=G\circ\varphi$. This means that $\varphi$ is a isomorphism if and only if $\varphi\circ H\circ\varphi^{-1}=G$. $\Box$

\begin{co}\label{Zc1} Two symmetric $m$-ary algebras of dimension $n$ are isomorphic if and only if their $n$-tuples are linearly conjugate.
\end{co}

Denote by
\bes
\mathrm{St}(H)=\{\a\in GL_n(k) \ | \ \a H \a^{-1}=H \}
\ees
the stabilizer of the homogeneous endomorphism $H$ in $GL_n(k)$.

\begin{lm}\label{Zl3}  The group of automorphisms $\mathrm{Aut}(P_H)$ of the polarization algebra $P_H$ is isomorphic to the stabilizer $\mathrm{St}(H)$ of $H$ in $GL_n(k)$.
\end{lm}
\Proof As in the proof of Lemma \ref{Zl2}, $\varphi\in \mathrm{Hom}(V,V)$ is an endomorphism of $P_H$ if and only if
 $\varphi(H(v))=H(\varphi(v))$ for all $v\in V$, i.e., $\varphi\circ H=H\circ\varphi$. $\Box$

Let $A$ be an arbitrary algebra with an $m$-ary $k$-linear operation $\langle \cdot,\cdot,\ldots,\cdot\rangle$.
Put $A^{(0)}=A$ and $A^{(i+1)}=\langle A^{(i)},A^{(i)},\ldots,A^{(i)}\rangle$ for all $i\geq 0$. The algebra $A$ is called {\em solvable} if there exists $p\geq 0$ such that  $A^{(p)}=0$.

\begin{lm}\label{Zl4} Let $H$ be a homogeneous of degree $m\geq 2$ endomorphism of the polynomial algebra $P_n$. Then $P_H^{(1)}=P_H$ if and only if $h_1,h_2,\ldots,h_n$ are linearly independent over $k$.
\end{lm}
\Proof Suppose that $P_H^{(1)}\neq P_H$ and let $f_1,f_2,\ldots,f_n$ be a basis for $P_H$ such that $f_s,f_{s+1},\ldots,f_n$ is a basis for $P_H^{(1)}$ where $s\geq 2$. Put
\bes
Y=x_1\otimes f_1+x_2\otimes f_2+\ldots+x_n\otimes f_n\in P_n\otimes P_H.
\ees
 Then
\bes
\langle Y,Y,\ldots,Y\rangle=g_1\otimes f_1+g_2\otimes f_2+\ldots+ g_n\otimes f_n,
\ees
where $G=(g_1,g_2,\ldots,g_n)$ is a homogeneous $n$-tuple of degree $m$. Notice that $g_1=\ldots g_{s-1}=0$ since
$f_s,f_{s+1},\ldots,f_n$ is a basis for $P_H^{(1)}$. By Lemma \ref{Zl2}, there exists $\varphi\in GL_n(k)$ such that
$\varphi\circ H=G\circ \varphi$. Comparing the first components, we get
\bes
\lambda_1h_1+\lambda_2h_2+\ldots+\lambda_nh_n=0,
\ees
if $\varphi(x_1)=\lambda_1x_1+\lambda_2x_2+\ldots+\lambda_nx_n$. Notice that $(\lambda_1,\lambda_2,\ldots,\lambda_n)\neq 0$ since $\varphi$ is invertible.

If $h_1,h_2,\ldots,h_n$ are linearly dependent over $k$ then it is not difficult to find $\varphi\in GL_n(k)$ such that the first component of $G=\varphi\circ H \circ\varphi^{-1}$ is zero. This means $P_G^{(1)}\neq P_G=P_H$. $\Box$

Let $A$ be an arbitrary algebra with an $m$-ary $k$-linear operation $\langle \cdot,\cdot,\ldots,\cdot\rangle$. Put  $A^1=A$, $A^m=\langle A,A,\ldots,A\rangle$, and
\bes
A^q=\sum_{i_1+i_2+\ldots+i_m=q} \langle A^{i_1},A^{i_2},\ldots,A^{i_m}\rangle
\ees
if there exists positive integers $i_1,i_2,\ldots,i_m$ such that $i_1+i_2+\ldots+i_m=q$ and $A^{i_1},A^{i_2},\ldots,A^{i_m}$ are defined.
Notice that $A^q$ is not defined for many $q$ if $m\geq 3$ (for example $A^2$ is not defined). The algebra $A$ is called {\em nilpotent} if there exists a number $t$ such that $A^q=0$ for all $q\geq t$.

This is a traditional definition of nilpotent algebras corresponding to the lower central series in group theory. We prefer a definition corresponding to the upper central series. Recall that an element $a\in A$ is called an {\em annihilator} of $A$ if
\bes
\langle b_1,\ldots,b_{i-1},a,b_{i+1},\ldots,b_m\rangle =0
\ees
for all $1\leq i\leq m$ and for all $b_1,\ldots,b_{i-1},b_{i+1},\ldots,b_m\in A$. Denote by $\mathrm{Ann}(A)$ the set of all annihilators of $A$.

Put $Z_0(A)=0$ and $Z_1(A)=\mathrm{Ann}(A)$. Suppose that $Z_i(A)$, where $i\geq 0$, is already defined and is an ideal of $A$. Then $Z_{i+1}(A)$ is defined as the full preimage of $\mathrm{Ann}(A/Z_i(A))$ under the natural homomorphism $A\rightarrow A/Z_i(A)$. The algebra $A$ is called {\em nilpotent} if there exists a natural $p$ such that $Z_p(A)=A$. It is well known that the two definitions of nilpotent algebras given above are equivalent.

Every nilpotent algebra is solvable but the converse is not true in general \cite{KBKA}.

An endomorphism $H=(h_1,h_2,\ldots,h_n)$ is called {\em strongly triangular} if $h_i\in k[x_1,\ldots,x_{i-1}]$ for all $i$.
An endomorphism $H=(h_1,h_2,\ldots,h_n)$ is called {\em strongly triangulable} if there exists $\a\in \mathrm{Aut}(\mathbb{A}^n)$ such that $\a H \a^{-1}$ is strongly triangular. A homogeneous endomorphism $H=(h_1,h_2,\ldots,h_n)$ is called {\em strongly linearly triangulable} if there exists $\a\in GL_n(k)$ such that $\a H \a^{-1}$ is strongly triangular.

\begin{lm}\label{Zl5} A polarization algebra $P_H$ is nilpotent if and only if $H$ is strongly linearly triangulable.
\end{lm}
\Proof Suppose that $A=P_H$ is nilpotent and let $p$ be the minimal positive integer such that $Z_p(A)=A$. Put
$k_i=\dim\,(A/Z_i(A)+1$ for all $1\leq i\leq p$.
Then  $1=k_p<k_{p-1}<\ldots <k_1\leq n$ and there exists a basis $(f)=\{f_1,f_2,\ldots,f_n\}$ for $A$ such that
$\{f_{k_i},\ldots,f_n\}$ is a basis for $Z_i(A)$ for all $1\leq i\leq p$. Let $G=(g_1,g_2,\ldots,g_n)$ be the $n$-tuple of $A$ corresponding to $(f)$. By induction on $n$ we prove that $G$ is strongly triangular. Notice that $G$ does not depend on $x_{k_1},\ldots,x_n$ since $\{f_{k_1},\ldots,f_n\}$ is a basis for $I=Z_1(A)=\mathrm{Ann}(A)$. Moreover, $G'=(g_1,g_2,\ldots,g_{k_1-1})$ is the $k_1-1$-tuple of $A/I$ corresponding to the basis $\{f_1+I,f_2+I,\ldots,f_{k_1-1}+I\}$. Besides, this basis satisfies the same conditions as $(f)$ with the numbers $1=k_p<k_{p-1}<\ldots <k_2\leq k_1-1$ for $A/I$. We may assume that $G'$ is strongly triangular since $k_1-1<n$. Consequently, $G$ is also strongly triangular.

Suppose that $H$ is strongly linearly triangulable. By Lemma \ref{Zl2}, we may assume that $H$ is strongly triangular. Then $H$ does not depend on $x_n$. Let $k_1$ be the minimal positive number such that $H$ does not depend on $x_{k_1},\ldots,x_n$. Denote by $I$ the span of $e_{k_1},\ldots,e_n$, where  $e_1,e_2,\ldots,e_n$ is the standard basis for $V=k^n$. It is easy to see that $I\subseteq \mathrm{Ann}(A)$. Moreover, $H'=(h_1,h_2,\ldots,h_{k_1-1})$ is the $k_1-1$-tuple of $A/I$ corresponding to the basis $\{e_1+I,e_2+I,\ldots,e_{k_1-1}+I\}$. Leading an induction on the dimension of $A$, we may assume that $A/I$ is nilpotent. Then $A$ is also nilpotent since $I\subseteq \mathrm{Ann}(A)$. $\Box$

\section{An approach to the Jacobian Conjecture}

\hspace*{\parindent}

In this section we describe an approach to the Jacobian Conjecture
 by L.M. Dru\.zkowski and K. Rusek \cite{DR85}, G. Gorni and G. Zampieri \cite{GZ96}, and by A.V. Yagzhev \cite{Yagzhev00-1}.

 It is well known \cite{Essen} that a homogeneous $n$-tuple $H$ has nilpotent Jacobian matrix
 \bes
 J(H)=(\partial_j(h_i))_{1\leq i,j\leq n}
 \ees
  if and only if the Jacobian
\bes
\mathrm{jac}(X-H)=\det(I-J(H))
\ees
 of the endomorphism $X-H=(x_1-h_1,x_2-h_2,\ldots,x_n-h_n)$ is $1$, where $I$ is the identity matrix of order $n$.
If $m\geq 2$ and $J(H)$ is nilpotent then the polynomial endomorphism $X-H=(x_1-h_1,x_2-h_2,\ldots,x_n-h_n)$ of $V=k^n$ is often called  a homogeneous Keller maps of degree $m$. Recall that the homogeneous Keller maps of degree $2$ are automorphisms \cite{Wang}. The well known Jacobian Conjecture is reduced to the study of homogeneous Keller maps of degree $3$ by Yagzhev \cite{Yagzhev} and Bass-Connel-Wright \cite{BCW}. Continuing this approach, L.M. Dru\.zkowski and K. Rusek \cite{DR85} first considered $H=H(X,X,X)$ as a trilinear operation. Later this algebra was called a ternary symmetric algebra by A.V. Yagzhev \cite{Yagzhev00-1}.

Let $k[V]$ be the algebra of all polynomial functions on the space $V=k^n$ and let $V^*$ be the space of all linear functions on $V$. If $e_1,e_2,\ldots,e_n$ is the standard basis for $V$ then $x_1,x_2,\ldots,x_n$ is the basis for $V^*$ dual to $e_1,e_2,\ldots,e_n$ with respect to the evaluation pairing
\bes
\langle u^*,u\rangle= u^*(u), \ \ \ u^*\in V^*, u\in V.
\ees
We have $k[V]=k[x_1,x_2,\ldots,x_n]=P_n$.
For every polynomial mapping $F : V\rightarrow V$  we define
$F^* : P_n\rightarrow P_n$ by $F^*(p)(v)=p(F(v))$ for all $p\in P_n$ and $v\in V$. It is well known that
$F^*$ is an automorphism of $P_n$ and if $F=(f_1,f_2,\ldots,f_n)$ then $F^*(x_i)=f_i$ for all $i$.

Notice that $X=(x_1,x_2,\ldots,x_n)=x_1\otimes e_1+x_2\otimes e_2+\ldots+x_n\otimes e_n$ is a generic element of the polarization algebra $P_H$. Then $X=(x_1,x_2,\ldots,x_n)$ generates the free algebra $P_H\langle X\rangle$ of the variety of algebras generated by $P_H$ in one free variable $X$. This section is mainly focused on the study of the free polarization algebra $P_H\langle X\rangle$.

\begin{lm}\label{Zl6} \cite{Yagzhev00-1} A polynomial endomorphism $X-H$ of the space $V=k^n$ is an automorphism if and only if the endomorphism of the free polarization algebra $P_H\langle X\rangle$ defined by $X\mapsto X-\langle X,X,\ldots,X\rangle$ is an automorphism.
\end{lm}
\Proof Put $F=X-H$ and consider the endomorphism $F^*$ of the polynomial algebra $P_n$.
Denote by $\Phi$ a unique $k$-linear extension of $F^*$ to the space
$k[X]\otimes P_H=k[X]^n$. Obviously, $\Phi$ is an endomorphism of the symmetric algebra since the product in $P_H$ does not depend on $x_1,\ldots,x_n$ and $\Phi$ keeps fixed every element of $k$. Moreover,
\bes
\Phi(X)=X-H=X-\langle X,X,\ldots, X\rangle.
\ees
Consequently, $\Phi$ induces an endomorphism $T$ of the free algebra $P_H\langle X\rangle$ defined by $X\mapsto X-\langle X,X,\ldots,X\rangle$. Let $T^{-1}$ be the formal inverse to $T$. Since $P_H\langle X\rangle$ is a homogeneous algebra, the existence of the formal inverse $T^{-1}$ to $T$ can be checked as in the case of polynomials \cite{Essen}. We have $T\circ T^{-1}=\mathrm{id}$, i.e., $T(T^{-1}(X))=X$. Notice that $T^{-1}(X)=(g_1,g_2,\ldots,g_n)=G$ where $g_i\in k[[X]]$ for all $i$. We get $F^*(G(X))=X$ since $T$ is induced by $\Phi$. This means that $F^*\circ G^*=\mathrm{id}$ and the formal inverses to $F^*$ and $T$ are the same. Consequently, $F^*$ is an automorphism if and only if $T$ is an automorphism. $\Box$

There are many possible ways to define an Engel condition for $m$-ary algebras if $m>2$. We choose the weakest one.
Let $A$ be an $m$-ary symmetric algebra. For any $a_1,a_2,\ldots,a_{m-1}\in A$ the multiplication operator
\bes
M(a_1,a_2,\ldots,a_{m-1}) : A \longrightarrow A
\ees
is defined by $x\mapsto \langle a_1,a_2,\ldots,a_{m-1},x\rangle$ for all $x\in A$. We say that $A$ is
an {\em Engel} algebra if there exists a positive integer $p$ such that
\bes
M(a,a,\ldots,a)^p=0
\ees
for all $a\in A$. If $m=2$ then this is the standard definition of a commutative Engel algebra.

The statement of the next lemma has been first noticed by G. Gorni and G. Zampieri \cite{GZ96} and proved by A.V. Yagzhev \cite{Yagzhev00-1} using the structural constants.
\begin{lm}\label{Zl7} \cite{GZ96,Yagzhev00-1}
 A polarization algebra $P_H$ is an Engel algebra if and only if the Jacobian matrix $J(H)$ is nilpotent.
\end{lm}
\Proof Consider two generic elements $X=(x_1,x_2,\ldots,x_n)$ and $Y=(y_1,y_2,\ldots,y_n)$ of $P_H$.
Let $d$ be the derivation of the polynomial algebra $k[X,Y]=k[x_1,\ldots,x_n,y_1,\ldots,y_n]$ such that $d(x_i)=y_i$ and $d(y_i)=0$ for all $i$. Denote by $D$ a unique $k$-linear extension of $d$ to the space
$k[X,Y]\otimes P_H=k[X,Y]^n$. Obviously, $D$ is a derivation of the symmetric algebra since the product in $P_H$ does not depend on $x_1,\ldots,x_n,y_1,\ldots,y_n$. Moreover, $D(X)=Y$ and $D(Y)=0$ and $P_H\langle X,Y\rangle$ is a $D$-closed free algebra in two variables. We have
\bes
\langle X,X,\ldots, X\rangle=H.
\ees
Applying $D$, we get
\bes
m M(X,X,\ldots,X)Y=D(H)=J(H)Y.
\ees
Therefore
\bes
m^p M(X,X,\ldots,X)^pY=J(H)^pY
\ees
for any positive integer $p$. Consequently, $J(H)$ is nilpotent if and only if $M(X,X,\ldots,X)$ is nilpotent. $\Box$

We give the definition of {\em nonassociative} (nonempty) words in one variable $X$ with respect to $m$-ary operation $\langle\cdot,\cdot,\ldots,\cdot\rangle$ following the definition for binary algebras  \cite{KBKA}:

(i) $X$ is a unique nonassociative word of length $1$;

(ii) If $w_1,w_2,\ldots,w_m$ are nonassociative words of length $k_1,k_2,\ldots,k_m\geq 1$, respectively, then
$\langle w_1,w_2,\ldots,w_m\rangle$ is a nonassociative word of length $k_1+k_2+\ldots+k_m$.

As in the binary case \cite{KBKA}, it is easy to show that every nonassociative word $w$ of length $l(w)>1$ has a unique decomposition $w=\langle w_1,w_2,\ldots,w_m\rangle$.

Notice that the free polarization algebra $P_H\langle X\rangle$ has grading
\bes
P_H\langle X\rangle= G_1\oplus G_2\oplus G_3\oplus\ldots \oplus G_s\oplus \ldots,
\ees
where $G_s$ the span of all nonassociative words in $X$ of length $(s-1)m-s+2$.

The statement $(i)$ of the following lemma was proved in \cite{DR85,Yagzhev00-1} and the statement $(ii)$ was proved in \cite{Yagzhev00-1}. The authors of \cite{GZ96} have used the word "nestings" instead of nonassociative words.

\begin{lm}\label{Zl8}
 Let $T$ be the endomorphism  of the free polarization algebra $P_H\langle X\rangle$ defined by $X\mapsto X-\langle X,X,\ldots,X\rangle$ and let $T^{-1}$ be the formal inverse to $T$. Then the following statements are true.

 $(i)$ $T^{-1}(X)= T_1+T_2+\ldots+T_s+\ldots$, where $T_1=X$ and
\bes
T_s=\sum_{i_1+\ldots+i_m=s} \langle T_{i_1},\ldots, T_{i_m}\rangle
\ees for all $s>1$.

 $(ii)$  $T_s$ is the sum off all nonassociative words in $X$ of degree $(s-1)m-s+2$ for all $s\geq 1$.
\end{lm}
\Proof Since $T^{-1}$ is the formal inverse to $T$ it follows that
\bes
\sum_{s\geq 1}T_s-\langle \sum_{s\geq 1}T_s,\ldots, \sum_{s\geq 1}T_s\rangle=X.
\ees
Comparing homogeneous parts, we get $(i)$.
Leading an induction on $s$, we can easily get $(ii)$ from $(i)$. $\Box$

Of course, the polynomials $T_s$ are different for different values of $m$.

\begin{co}\label{Zc2} Let $H$ be a homogeneous $n$-tuple of elements of $P_n$ of degree $m\geq 2$.
Then the polynomial mapping $X-H$ of $V=k^n$ is an automorphism if and only if there exists a positive integer $p=p(H)$ such that the polarization algebra $P_H$ satisfies all identities $T_s=0$ where $s\geq p$.
\end{co}
\Proof By Lemma \ref{Zl6}, $X-H$ is an automorphism if and only if the endomorphism $T$ of the free polarization algebra $P_H\langle X\rangle$ defined by $X\mapsto X-\langle X,X,\ldots,X\rangle$ is an automorphism. By Lemma \ref{Zl8}, $T$ is automorphism if and only if there exists $p$ such that $T_s=0$ for all $s\geq p$. This means that $P_H$ satisfies the identities $T_s=0$ for all $s\geq p$. $\Box$
\begin{co}\label{Zc3} The Jacobian Conjecture for polynomial algebras is true if and only if for any finite dimensional symmetric $m$-ary Engel algebra $A$ there exists a positive integer $p=p(A)$ such that  $A$ satisfies all identities $T_s=0$ where $s\geq p$.
\end{co}
\Proof The Jacobian Conjecture is true if and only every polynomial endomorphism $X-H$ with nilpotent Jacobian matrix $J(H)$ is an automorphism. By Lemma \ref{Zl1}, every symmetric algebra $A$ is isomorphic to $P_H$. By Lemma \ref{Zl7},  $A$ is an Engel algebra if and only if $J(H)$ is nilpotent. Corollary \ref{Zc2} finishes the proof. $\Box$

Combining Corollary \ref{Zc2} and the Jacobian Conjecture for quadratic Keller maps \cite{Wang} we get the next statement.
\begin{co}\label{Zc4} For any finite dimensional commutative Engel algebra $A$ there exists a positive integer $p=p(A)$ such that  $A$ satisfies all identities $T_s=0$ where $s\geq p$.
\end{co}

\section{Albert's problem and the homogeneous dependence problem}

\hspace*{\parindent}

Recall that an algebra $A$ (not necessarily associative) over a field $k$ is called  {\em power-associative} \cite{KBKA}  if each element of $A$ generates an
associative subalgebra. For power-associative algebras, an element $a$ is called {\em nilpotent} if $a^n=0$ for some $n$ and an algebra is called {\em nil} if every its element is nilpotent.

Albert's conjecture was \cite{Albert} that every commutative power-associative finite dimensional nilalgebra is nilpotent. Gerstenhaber and Myung \cite{GM} proved that every commutative power-associative nil algebra  of dimension $\leq 4$ is nilpotent. In \cite{Suttles} Suttles gave an example of a commutative power-associative algebra of dimension 5 which is not nilpotent but solvable. As we mentioned in the introduction, the automorphism (\ref{Zf1}) corresponds to the Suttles' example. Notice that this automorphism can be written in the form $X-H$ where $X=(x,y,z,s,t)$ and $H=(ys,xt-zs,yt,0,0)$. Then $J(H)$ is nilpotent since $X-H$ is an automorphism. Consequently, $P_H$ is a binary commutative $5$ dimensional algebra with an Engel identity. By the definition of the product in $P_H$ we have
$X\cdot X=H$. Consequently,
\bes
X\cdot X=(xe_1+ye_2+ze_3+se_4+te_5)\cdot (xe_1+ye_2+ze_3+se_4+te_5)\\
=x^2e_1\cdot e_1+ 2xye_1\cdot e_2+2xze_1\cdot e_3+2xse_1\cdot e_4+2xte_1\cdot e_5\\
+y^2e_2\cdot e_2+ 2yze_2\cdot e_3+2yse_2\cdot e_4+2yte_2\cdot e_5
+z^2e_3\cdot e_3\\
+2zse_3\cdot e_4+2zte_3\cdot e_5 + s^2e_4\cdot e_4+2ste_4\cdot e_5+t^2e_5\cdot e_5\\
=H=yse_1+(xt-zs)e_2+yte_3.
\ees
Comparing coefficients in monomials, we get the following nonzero products
\bes
e_2\cdot e_4=1/2e_1, e_1\cdot e_5=1/2e_2, e_3\cdot e_4=-1/2e_2, e_2\cdot e_5=1/2e_3.
\ees
Using this, we can easily get that $X^2\cdot X^2=0$ and $(X^2\cdot X)\cdot X=0$. It is well known (see, for example \cite{Fernandez04}) that a commutative algebra
over a field of characteristic zero is power-associative if and only if the identity
\bes
X^2X^2=(X^2X)X
\ees
holds. Consequently, the polarization algebra $P_H$ is power-associative. It is easy to check that $P_H$ is not nilpotent but solvable.

Now consider the automorphism (\ref{Zf2}). This automorphism can be written in the form $X-H$ where $X=(x,y,s,t)$ and $H=((ys-xt)s,(ys-xt)t,0,0)$. Then $P_H$ is a ternary symmetric $4$ dimensional algebra with an Engel identity. We have
\bes
\langle X,X,X\rangle=x^3\langle e_1,e_1,e_1\rangle+3x^2y\langle e_1,e_1,e_2\rangle+3x^2s\langle e_1,e_1,e_3\rangle
+3x^2t\langle e_1,e_1,e_4\rangle\\
+3xy^2\langle e_1,e_2,e_2\rangle+6xys\langle e_1,e_2,e_3\rangle+6xyt\langle e_1,e_2,e_4\rangle
+3xs^2\langle e_1,e_3,e_3\rangle\\
+6xst\langle e_1,e_3,e_4\rangle+3xt^2\langle e_1,e_4,e_4\rangle
+y^3\langle e_2,e_2,e_2\rangle+3y^2s\langle e_2,e_2,e_3\rangle\\
+3y^2t\langle e_2,e_2,e_4\rangle
+3ys^2\langle e_2,e_3,e_3\rangle
+6yst\langle e_2,e_3,e_4\rangle
+s^3\langle e_3,e_3,e_3\rangle\\
+3s^2t\langle e_3,e_3,e_4\rangle
+3st^2\langle e_3,e_4,e_4\rangle
+3t^3\langle e_4,e_4,e_4\rangle\\
=H=(ys-xt)se_1+(ys-xt)te_2.
\ees
From this we get the following nonzero products
\bes
\langle e_2,e_3,e_3\rangle=1/3 e_1, \langle e_1,e_3,e_4\rangle=-1/6 e_1,\\
 \langle e_2,e_3,e_4\rangle=1/6 e_2, \langle e_1,e_4,e_4\rangle=-1/3 e_2.
\ees
Using this, it is easy to check that the free polarization algebra $P_H\langle X\rangle$ is a $2$ dimensional nilpotent algebra with a linear basis $X,\langle X,X,X\rangle$, Consequently, $P_H\langle X\rangle$ is associative.
This means that $P_H$ is power-associative. It is easy to check that $P_H$ is not nilpotent but solvable.

A modified version of the Albert's problem  was formulated in
 \cite[Problem 1.1]{Dnestr}. \\
{\bf Albert's problem I:}
Is every finite dimensional commutative power-associative nilalgebra solvable?

This problem is still open.  It has been proved \cite{CHP,CJ} that commutative power-associative nil algebras  of dimension $\leq 8$ are solvable. This problem is also equivalently related to possessing of symmetric associative bilinear forms \cite{Arenas}. Commutative power-associative algebras with a nil basis were studied in \cite{Dedkov,Shestakov}.

An algebra $A$ is called a {\em nilalgebra of bounded nilindex} if there exists a positive integer $t$
such that $x^
t = 0$  for all $x\in A$ and all possible distribution of parenthesis.
Many authors started to study the following version of the Albert's problem.\\
{\bf Albert's problem II}: Is every finite dimensional commutative nilalgebra (not necessarily power-associative) solvable?

The latest results say that all commutative nilalgebras  of dimension $\leq 7$ are solvable \cite{Fernandez04,Fernandez10}.

 In 1960 M. Gerstenhaber proved that \cite{Gerstenhaber60} every commutative nilalgebra of bounded index is an Engel algebra.
 Many authors also started to study the following version of the Albert's problem.\\
{\bf Albert's problem III}: Is every finite dimensional commutative Engel algebra solvable?

It is proved that every finitely generated commutative algebra with the Engel identity $x(x(xy))=0$ is solvable \cite{CH} and every finite dimensional commutative algebra
with the same identity is nilpotent \cite{Fernandez09}. The classification of homogeneous quadratic automorphisms in dimension $5$ given in \cite{Sun} can be considered as a classification of $5$ dimensional commutative Engel algebras since all quadratic Keller maps are automorphisms \cite{Wang}.

Of course, all three versions of the Albert's problem can be considered for $m$-ary symmetric algebras. Let $A$ be an arbitrary $m$-ary algebra and $A\langle X\rangle$ be the free algebra in one free variable $X$ of the variety of algebras generated by $A$. The algebra $A$ is called {\em power-associative} if $A\langle X\rangle$ is associative and is called {\em nil of bounded index} if $A\langle X\rangle$ is nilpotent.

The next lemma is an analogue of the result by Gerstenhaber \cite{Gerstenhaber60} mentioned above .
\begin{lm}\label{Zl9}
Let $A$ be an arbitrary $m$-ary symmetric nilalgebra of bounded index. Then $A$ is an Engel algebra.
\end{lm}
\Proof Let $A\langle X\rangle$ be the free algebra in one free variable $X$ of the variety of algebras generated by $A$.
Consider the endomorphism $T$ of $A\langle X\rangle$ defined by $X\mapsto X-\langle X,X,\ldots,X\rangle$. Then $T$ is an automorphism since $A\langle X\rangle$ is nilpotent. Therefore the Jacobian matrix $J(T)=1-mM(X,\ldots,X)$ of $T$ is invertible \cite{Um11}. Consequently (see, for example \cite{Essen}), $M(X,\ldots,X)$ is nilpotent. $\Box$

Now recall the formulation of the homogeneous dependence problem. \\
{\bf HDP(m)}. Let $H=(h_1,h_2,\ldots,h_n)$ be a homogeneous $n$-tuple of elements of degree $m$ of the polynomial algebra $P_n$.
If $J(H)$ is nilpotent then is it true that $h_1,h_2,\ldots,h_n$ are linearly dependent over $k$?

\begin{theor}\label{Zt1} Let $m\geq 2$ be an arbitrary integer.
The homogeneous dependence problem  HDP(m) is true if and only if every finite dimensional $m$-ary symmetric Engel algebra is solvable.
\end{theor}
\Proof Suppose that every finite dimensional $m$-ary symmetric Engel algebra is solvable.
Let $H=(h_1,h_2,\ldots,h_n)$ be a homogeneous $n$-tuple of elements of degree $m$ of $P_n$ with nilpotent Jacobian matrix $J(H)$. Then $P_H$ is a  finite dimensional $m$-ary symmetric Engel algebra by Lemma \ref{Zl7}. Consequently, $P_H$ is solvable. Lemma  \ref{Zl4} gives that $h_1,h_2,\ldots,h_n$ are linearly dependent over $k$.

Suppose that the homogeneous dependence problem  HDP(m) is true. Let $A$ be an arbitrary finite dimensional $m$-ary symmetric Engel algebra and $H=(h_1,h_2,\ldots,h_n)$ be an $n$-tuple of $A$. Then $A$ is isomorphic to $P_H$ by Lemma  \ref{Zl1}. By Lemma \ref{Zl7}, $J(H)$ is nilpotent. Then $h_1,h_2,\ldots,h_n$ are linearly dependent since HDP(m) is true. Lemma  \ref{Zl4} gives that $A^{(1)}\neq A$. Consequently, $\mathrm{dim}\,A^{(1)}< \mathrm{dim}\,A$. Leading an induction on the dimension of $A$, we may assume that $A^{(1)}$ is solvable. Then $A$ is also solvable. $\Box$

Let's consider some results on the homogeneous dependence problem. Some positive results in small dimensions can be found in \cite{Essen}. The homogeneous dependence problem HDP(m) is negatively solved for all $m\geq 3$ by M. de Bondt \cite{Bondt06,Bondt09}. A counterexample of dimension $9$ for $m=3$ is given in \cite{Zampieri08}.
\begin{co}\label{Zc5}
An analogue of the Albert's problem III for $m$-ary algebras has a negative solution  if $m\geq 3$.
\end{co}
The real algebraic case $m=2$ remains open. The methods of construction of counter-examples for $m\geq 3$ might be useful in binary case. In fact, M. de Bondt \cite{Bondt06} started his construction from the automorphism (\ref{Zf2}). It might be interesting try to get a binary counter-example starting from the binary automorphism (\ref{Zf1}).

Of course, specialists in affine algebraic geometry did not consider the power-associativity. Occasionally, $5$-ary algebras of dimension $n\geq 6$ corresponding to the automorphisms given in Theorem 2.1 in \cite{Bondt06} are power-associative and every one generated subalgebra of these algebras are nilpotent and at most $2$ dimensional.
\begin{co}\label{Zc6}
Analogues of Albert's problems I and II for $m$-ary algebras have negative solutions if $m\geq 5$.
\end{co}

\section{Comments and problems}

\hspace*{\parindent}

To prove a result on finite dimensional algebras we usually use an induction on the dimension of an algebra.
An advantage of considering $m$-ary algebras for all $m\geq 2$ is that we get one more induction to prove results on  finite dimensional symmetric nil algebras of bounded index. Let $A$ be an $n$ dimensional $m$-ary symmetric nil algebra of bounded index. Let $A\langle X\rangle$ be the free algebra in one free variable $X$ of the variety of algebras generated by $A$. Then $A\langle X\rangle$ is nilpotent and finite dimensional. Put
\bes
H(X)=\langle \langle X,X,\ldots,X\rangle,X,\ldots,X\rangle
\ees
and consider the free polarization algebra $P_H\langle X\rangle$. Obviously,  $P_H\langle X\rangle$ is associative if $A\langle X\rangle$ is associative and $P_H\langle X\rangle$ is nilpotent if $A\langle X\rangle$ is nilpotent.
Moreover,  $\dim\,P_H\langle X\rangle < \dim\,A\langle X\rangle$.  This allows to lead an induction on $\dim\,A\langle X\rangle$.

The critical case is when $\dim\,A\langle X\rangle=2$. In this case $H=0$ and all elements of $A\langle X\rangle$ except $X$ and $\langle X,X,\ldots,X\rangle$ are zeroes. Examples given in Theorem 2.1 in \cite{Bondt06} correspond exactly to this case. They do not allow to use an induction on $\dim\,A\langle X\rangle$ to solve Albert's problems I and II. But this induction may be useful in solving some other problems.

By Lemma \ref{Zl8}, the structure of the free polarization algebra $P_H\langle X\rangle$ is very important in the study of the Jacobian Conjecture. Notice that the polarization algebras corresponding to (\ref{Zf1}) and (\ref{Zf2}) are not nilpotent but the corresponding free polarization algebras in one variable are nilpotent.  An example of a homogeneous $4$-tuple $H$ of degree $3$ with nilpotent Jacobian matrix $J(H)$ such that $P_H\langle X\rangle$ is not nilpotent is given by G. Gorni and G. Zampieri \cite{GZ96}. In fact, the corresponding automorphism
\bes
(x+s(xt-ys),y+t(xt-ys),s+t^3,t)
\ees
is obtained by only one elementary transformation from (\ref{Zf2}). The same result we can get if we consider the automorphism
\bes
(x-ys,y+zs-xt,z-yt,s+t^2,t)
\ees
obtained from (\ref{Zf1}) by one elementary transformation. In this case $P_H$ becomes a binary algebra. So, elementary transformations hurt the nilpotency and the power-associativity.
In both cases $P_H\langle X\rangle$ is solvable since $P_H$ is solvable. It seems that elementary transformations do not hurt solvability of $P_H\langle X\rangle$.

All examples in \cite{Bondt06} based on quasi-translations and elementary transformations. Recall that $X-H$ is called a quasi-translation if $(X-H)^{-1}=X+H$. In this case $P_H\langle X\rangle$ is two dimensional nilpotent algebra.
\begin{prob}\label{prob1}
Is $P_H\langle X\rangle$ solvable if $J(H)$ is nilpotent?
\end{prob}

The weakest form of nil elements in nonassociative algebras can be defined in the following way. Let $A$ be an $m$-ary algebra and $a\in A$. Put $a^{(0)}=a$ and $a^{(i+1)}=\langle a^{(i)},a^{(i)},\ldots,a^{(i)}\rangle$. An element $a\in A$ is called weak nil if there exists a positive integer $p$ such that $a^{(p)}=0$. An algebra is called weak nil if every its element is weak nil.
\begin{prob}\label{prob2}
Is every finite dimensional symmetric $m$-ary Engel algebra is weak nil?
\end{prob}
This problem can be equivalently formulated in purely affine algebraic geometry language.
\begin{prob}\label{prob3}
Let $H$ be a arbitrary $n$-tuple with nilpotent Jacobian matrix $J(H)$ and $H(0)=0$. Then is the polynomial endomorphism $H$ nilpotent?
\end{prob}
Obviously, this problem can be reduced to homogeneous mappings.

\bigskip

\begin{center}
{\bf\large Acknowledgments}
\end{center}

\hspace*{\parindent}

I am grateful to Max-Planck Institute f\"ur Mathematik for their
hospitality and excellent working conditions, where part of this work has been done.

\end{document}